\documentclass[12pt,reqno]{article}

\usepackage[usenames]{color}
\usepackage{amssymb}
\usepackage{amsmath}
\usepackage{amsthm}
\usepackage{amsfonts}
\usepackage{amscd}
\usepackage{graphicx}
\usepackage[shortlabels]{enumitem}

\usepackage[colorlinks=true,
linkcolor=webgreen,
filecolor=webbrown,
citecolor=webgreen]{hyperref}

\definecolor{webgreen}{rgb}{0,.5,0}
\definecolor{webbrown}{rgb}{.6,0,0}

\usepackage{color}
\usepackage{fullpage}
\usepackage{float}

\usepackage{graphics}
\usepackage{latexsym}
\usepackage{epsf}
\usepackage{breakurl}

\def\Enn{\mathbb{N}}

\newcommand{\seqnum}[1]{\href{https://oeis.org/#1}{\rm \underline{#1}}}

\def\suchthat{\, : \, }

\begin{document}

\theoremstyle{plain}
\newtheorem{theorem}{Theorem}
\newtheorem{corollary}[theorem]{Corollary}
\newtheorem{lemma}[theorem]{Lemma}
\newtheorem{proposition}[theorem]{Proposition}

\theoremstyle{definition}
\newtheorem{definition}[theorem]{Definition}
\newtheorem{example}[theorem]{Example}
\newtheorem{conjecture}[theorem]{Conjecture}

\theoremstyle{remark}
\newtheorem{remark}[theorem]{Remark}

\def\kaw{Kawsumarng et al.}

\title{Sumsets of Wythoff Sequences, Fibonacci Representation, and Beyond}

\author{Jeffrey Shallit\\
School of Computer Science \\
University of Waterloo \\
Waterloo, ON  N2L 3G1 \\
Canada \\
\href{mailto:shallit@uwaterloo.ca}{\tt shallit@uwaterloo.ca} }

\maketitle

\begin{abstract}
Let $\alpha = (1+\sqrt{5})/2$ and define
the lower and upper Wythoff sequences by
$a_i = \lfloor i \alpha \rfloor$,
$b_i = \lfloor i \alpha^2 \rfloor$
for $i \geq 1$.
In a recent interesting paper, \kaw~proved a number of results
about numbers representable as sums of the form $a_i + a_j$,
$b_i + b_j$, $a_i + b_j$, and so forth.  In this paper I show
how to derive all of their results, using one
simple idea and existing free software called {\tt Walnut}.
The key idea is that for each of their sumsets, there is a
relatively small
automaton accepting the Fibonacci representation of the
numbers represented.
I also show how the automaton approach can easily prove
other results.
\end{abstract}

\section{Introduction}

Let $\alpha = (1+\sqrt{5})/2$ and define
the lower and upper Wythoff sequences by
$a_i = \lfloor i \alpha \rfloor$,
$b_i = \lfloor i \alpha^2 \rfloor$
for $i \geq 1$.
The first few terms of these sequences are as follows:
\begin{center}
\begin{tabular}{c|c|cccccccccc} 
& $i$ & 1 & 2 & 3 & 4 & 5 & 6 & 7 & 8 & 9 & 10 \\
\hline
lower Wythoff & $a_i$ & 1 & 3 & 4 & 6 & 8 & 9 & 11 & 12 & 14 &  16 \\
upper Wythoff & $b_i$ & 2 & 5 & 7 & 10 & 13 & 15 & 18 & 20 & 23 &  26
\end{tabular}
\end{center}
They are, respectively, sequences 
\seqnum{A000201} and \seqnum{A001950} in the
{\it On-Line Encyclopedia of Integer Sequences} (OEIS)
\cite{Sloane:2020}.

In a recent interesting paper, \kaw~\cite{Kaws}
proved a number of results in the additive number theory of the
Wythoff sequences; in particular,
about numbers representable as sums of the form $a_i + a_j$,
$b_i + b_j$, $a_i + b_j$, and so forth.  In this paper I show
how to derive all of their results, using one
simple idea and existing free software called {\tt Walnut}.
The key idea is that for each of their sumsets, there is a
relatively small
finite automaton accepting the Fibonacci representation of the
numbers represented, and one can read off the properties of
the represented numbers directly from it.
The idea of using automata to solve problems in additive number
theory has been explored previously in
\cite{Bell,Rajasekaran}.

\section{Fibonacci representation}

Let's start with Fibonacci representation (also called 
Zeckendorf representation).  The Fibonacci numbers
are given by $F_0 = 0$, $F_1 = 1$, and 
$F_i = F_{i-1} + F_{i-2}$ for $i \geq 2$.  The
Fibonacci representation
of a number $n$ is a representation
$$ n = \sum_{2 \leq i \leq t} e_i F_i ,$$
where $e_i \in \{0,1\}$ and we impose the conditions
that $e_t = 1$ and $e_i e_{i+1} = 0$ for $2 \leq i < t$.
As is well known, every positive integer has such
a representation, and this representation is unique
\cite{Lekkerkerker:1952,Zeckendorf:1972}.   
This representation of $n$ is usually abbreviated by
the bitstring $e_t e_{t-1} \cdots e_2$, and is written
as $(n)_F$.    Thus, for example,
$(43)_F = 10010001$, because $43 = F_9 + F_6 + F_2 = 34 + 8 + 1$.
Note that, conventionally, a Fibonacci representation is
written with the most significant bit at the left.
The integer $0$ is represented by the empty string $\epsilon$.

\section{Finite automata}

Next, let's talk about finite automata.   A finite automaton $A$
is a simple machine model having a finite number of states,
an initial state, transition rules that describe how inputs
cause the machine to move from state to state, and a notion of
``final'' or ``accepting states''.    An input is accepted if,
starting from the initial state, and following the transition
rules corresponding to the symbols of the input, the machine
reaches a final state.    For more details, see any introductory
textbook on automata theory, such as \cite{Hopcroft&Ullman:1979}.

By associating a number with the string (or word) given by its
Fibonacci representation, the strings accepted by an automaton $A$
can be considered as a subset $S= S(A)$ of the natural numbers $\Enn$.
When an automaton is interpreted in this way, we call it a
{\it Fibonacci automaton} and the corresponding set is called
{\it Fibonacci automatic}.  Many properties of Fibonacci-automatic
sets are known; for example, see
\cite{Mousavi&Schaeffer&Shallit:2016,
Du&Mousavi&Rowland&Schaeffer&Shallit:2017,
Du&Mousavi&Schaeffer&Shallit:2016}.

Let's look at an example.   Two famous classical identities on
Fibonacci numbers are as follows:
\begin{align*}
F_2 + F_4 + F_6 + \cdots + F_{2t} &= F_{2t+1} - 1; \\
F_3 + F_5 + \cdots + F_{2t-1}&= F_{2t} - 1;
\end{align*}
they hold for $t \geq 1$.
These identities, then, can be viewed as giving the Fibonacci
representation of numbers of the form $F_n - 1$ for $n \geq 2$;
that is, the set $ S = \{ 1,2,4,7,12, \ldots \}$.   The set
$S$ can then by represented by the following automaton:
\begin{figure}[H]
\begin{center}
\includegraphics[width=5in]{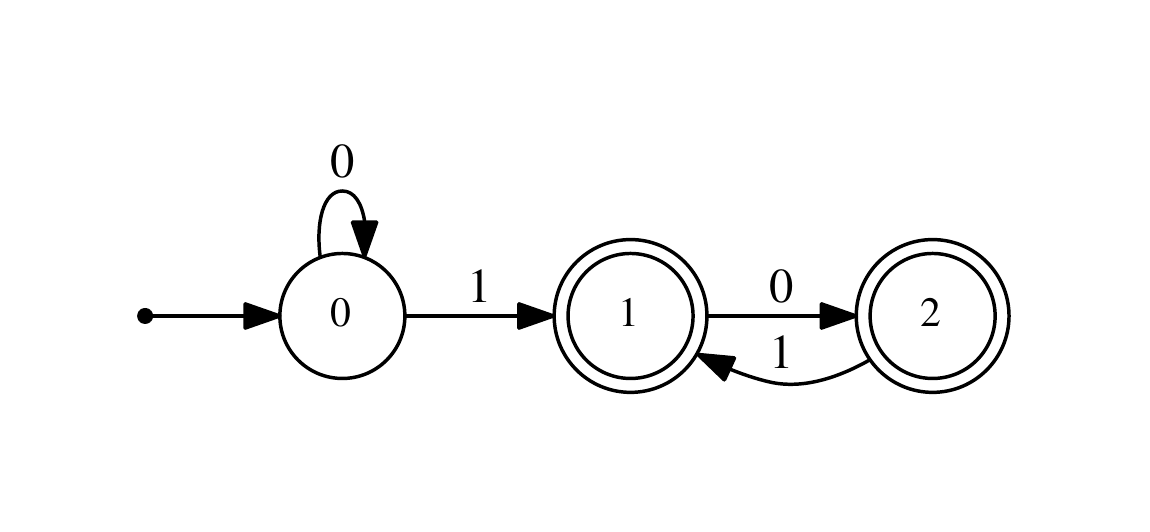}
\end{center}
\caption{Fibonacci automaton for the set $S = \{ 1,2,4,7,12, \ldots \}$}
\label{fig1}
\end{figure}
Here the initial state is denoted by headless arrow entering,
and the accepting states are denoted by double circles.  To use
this automaton, express $n$ in Fibonacci representation,
start at the initial state labeled $0$, and follow the arrows
according to the bits of $(n)_F$.   If the automaton ends up
in an accepting state, then $n \in S$.   If at some point
there is no possible move, or if the automaton ends in a
non-accepting state, then $n \not\in S$.

\section{Wythoff sequences}

The values of
the lower and upper Wythoff sequences have 
distinctive Fibonacci representations that
were determined by Silber \cite{Silber:1976}:
\begin{theorem}
\leavevmode
\begin{itemize}
\item $n$ belongs to the lower Wythoff sequence iff $(n-1)_F$
ends in $0$;
\item $n$ belongs to the upper Wythoff sequence iff $(n-1)_F$
ends in $01$.
\end{itemize}
For the purpose of this theorem, we pretend that every Fibonacci
representation is preceded by at least one $0$.
\end{theorem}

It follows from this (we will see exactly how below) that
there is a very simple
Fibonacci automaton for the set of lower Wythoff numbers,
and another for the set of upper Wythoff numbers.   We give
them below.  \kaw~\cite{Kaws} called these
sets $B(\alpha)$ and $B(\alpha^2)$, respectively, but we
shall call them $L$ and $U$ (for lower and upper).  
\begin{figure}[H]
\begin{center}
\includegraphics[width=3.1in]{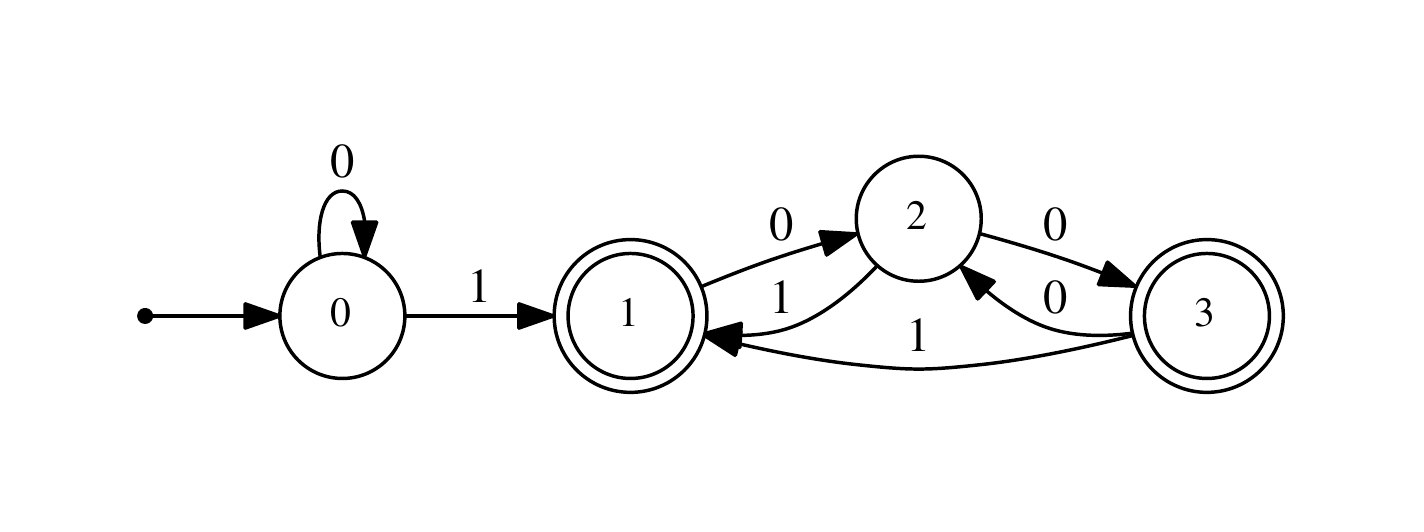}
\includegraphics[width=3.1in]{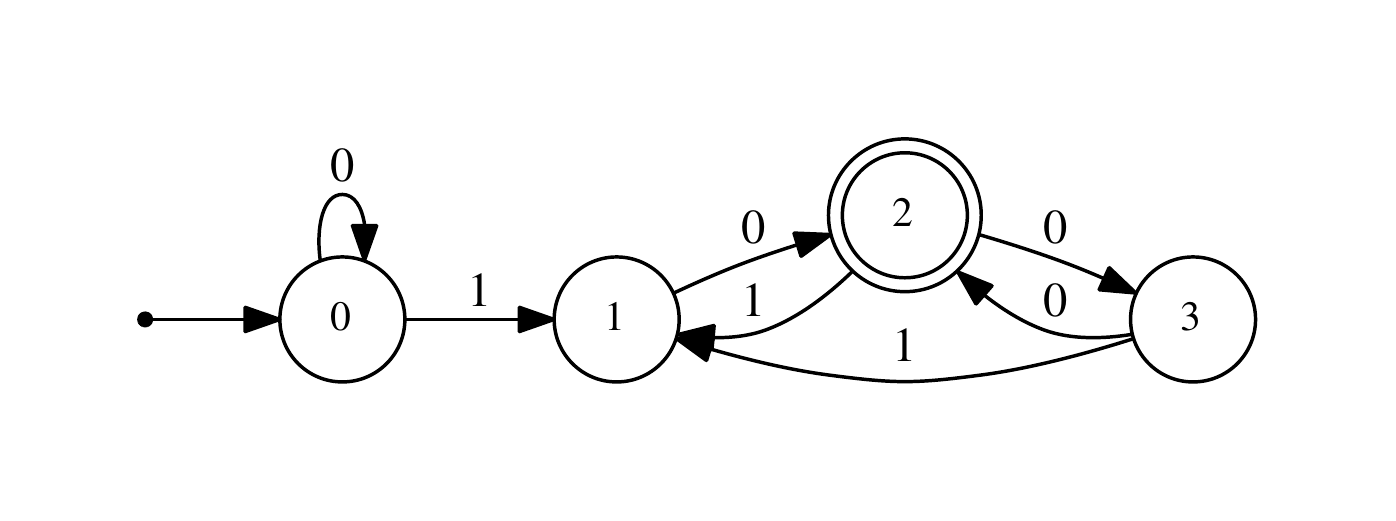}
\end{center}
\caption{Fibonacci automata for the sets $L$ (left) and $U$ (right)}
\label{fig2}
\end{figure}
Similarly, \kaw~defined $B_0 (\alpha) = B_(\alpha) \cup \{0\}$ and
$B_0 (\alpha^2) = B(\alpha^2) \cup \{0\}$, but we will call
these $L_0$ and $U_0$, respectively.

Recall that if $S, T \subseteq \Enn$, then
$S+T$ denotes the sumset $\{ s + t \suchthat s \in S,\ t \in T \}$.
Now that we have automata for $L$ and $U$, we can find automata
for their sumsets, using
the following theorem
\cite{Mousavi&Schaeffer&Shallit:2016}:
\begin{theorem}
Given a first-order formula $F$ involving a Fibonacci-automatic
set (or sets), addition, comparison, logical operators, and
the universal and existential quantifiers, there is an automaton
accepting the Fibonacci representations of those settings of
the free variables that make $F$ true.
\end{theorem}

Not only does
such an automaton exist; there is also an algorithm to
produce it.  This algorithm has been implemented in the
free software called {\tt Walnut} and written by Hamoon Mousavi
\cite{Mousavi:2016}.

For example, if we wish to compute the Fibonacci automaton for
the sumset $L + U$, we only need to write the first-order
formula
$$ \phi:  \exists a, b \ (n=a+b) \ \wedge \ a \in L \ \wedge\ b \in U, $$
translate it into {\tt Walnut}'s syntax,
\begin{verbatim}
def lplusu "?msd_fib E a,b (n=a+b) & $lower(a) & $upper(b)":
\end{verbatim}
and have {\tt Walnut} compute the corresponding 
automaton recognizing the language
$$ L(\phi) = \{ (n)_F \suchthat \phi \text{ holds } \} .$$
{\tt Walnut}'s syntax is more or less self-explanatory;
we note
that in {\tt Walnut}, the symbol
{\tt E} represents the existential quantifier and the
cryptic notation {\tt ?msd\_fib} indicates that we want to evaluate
the predicate using the Fibonacci representation of integers.
Here
\begin{itemize}
\item {\tt lower} is {\tt Walnut}'s name for $L := B(\alpha)$;
\item {\tt lower0} is {\tt Walnut}'s name for $L_0 := B_0 (\alpha) = B(\alpha) \cup \{0\} $;
\item {\tt upper} is {\tt Walnut}'s name for $U := B(\alpha^2)$;
\item {\tt upper0} is {\tt Walnut}'s name for $U_0 := B_0 (\alpha^2) = B(\alpha^2) \cup \{0\} $.
\end{itemize}

\section{Finding automata for the sets $L, U, L_0, U_0$}

As we have seen, $n \in L$ iff $(n-1)_F$ ends in $0$, and
$n \in U$ iff $(n-1)_F$ ends in $01$.   The easiest way to specify
these requirements is with a regular expression \cite{Hopcroft&Ullman:1979}.
In {\tt Walnut} we do this as follows:
\begin{verbatim}
reg end0 msd_fib "(0|1)*0":
reg end01 msd_fib "1|((0|1)*01)":
\end{verbatim}
The alert reader will notice that we have allowed an arbitrary
number of leading zeros at the beginning of our specifications.
This is needed in {\tt Walnut} for technical reasons.
The straight bar $|$ represents logical ``or''.

Now we can easily create
{\tt Walnut} predicates that assert that the Fibonacci
representation of $n-1$ matches the given expressions.
\begin{verbatim}
def lower "?msd_fib Em $end0(m) & n=m+1":
def lower0 "?msd_fib $lower(n) | (n=0)":
def upper "?msd_fib Em $end01(m) & n=m+1":
def upper0 "?msd_fib $upper(n) | (n=0)":
\end{verbatim}
When we run {\tt Walnut} on these predicates, we get the automata
given in Figure~\ref{fig2}, as well as two automata for $L_0$ and
$U_0$ (not depicted).

\section{Theorem 3.1 of \kaw}
We can now prove Theorem 3.1 of \kaw, namely:
\begin{theorem}
\leavevmode
\begin{enumerate}[(i)]
\item $L+L = \Enn - \{0,1,3\}$;
\item $L_0+L = \Enn - \{0\}$.
\end{enumerate}
\end{theorem}
(The slight difference in notation between their paper and ours is because
we use $\Enn = \{ 0,1,2,\ldots \}$, while
\kaw~used $\Enn = \{ 1,2,3, \ldots \}$.)

\begin{proof}
We use the {\tt Walnut} commands
\begin{verbatim}
def thm31i "?msd_fib ~Ea,b (n=a+b) & $lower(a) & $lower(b)":
def thm31ii "?msd_fib ~Ea,b (n=a+b) & $lower0(a) & $lower(b)":
\end{verbatim}
which compute automata for the {\it complements} 
$\Enn - (L+L)$ and $\Enn - (L_0 + L)$.
These commands produce the automata depicted below.  
\begin{figure}[H]
\begin{center}
\includegraphics[width=4in]{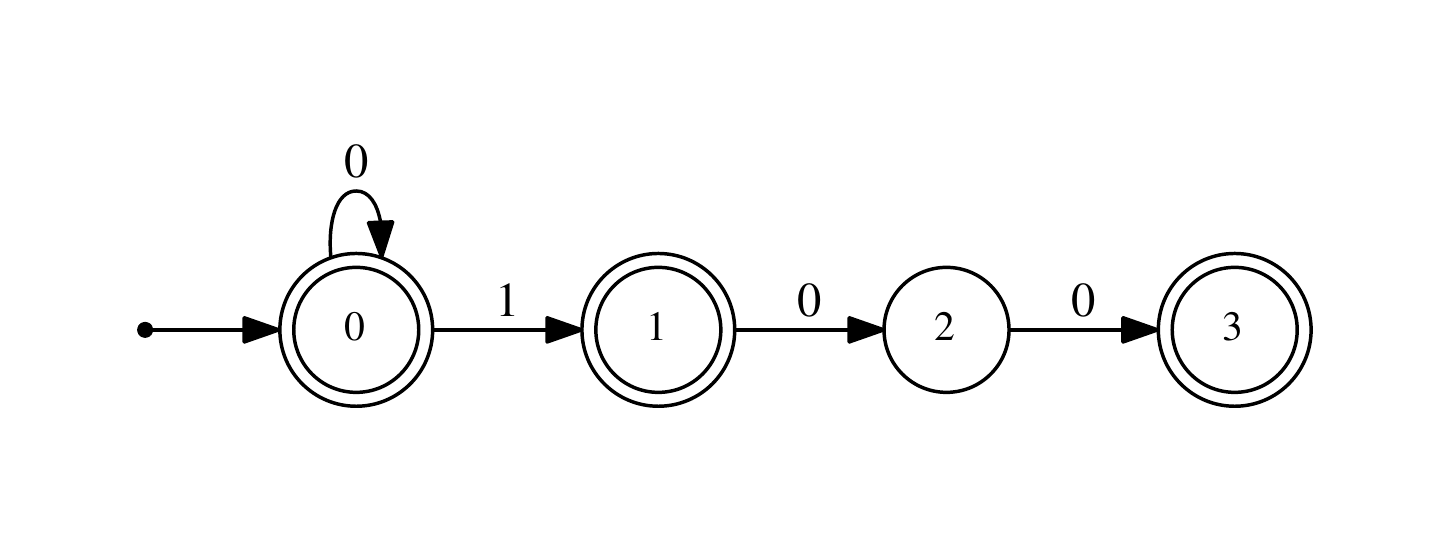}
\includegraphics[width=2in]{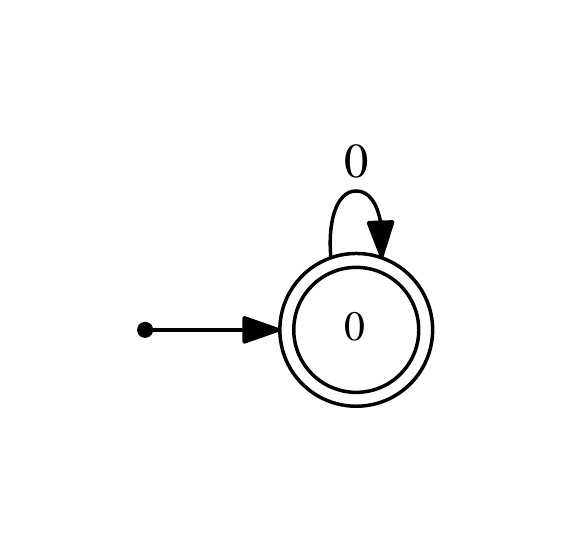}
\end{center}
\caption{Automata for $\Enn - (L+L)$ (left)
and $\Enn - (L_0 + L)$ (right)}
\end{figure}
\end{proof}

\begin{remark}
We point out that here the theorem-proving software {\tt Walnut} is not
just confirming something that could be conjecture from numerical
evidence; rather, inspection of the automata produced here
actually produces the {\it statement of the theorem itself\/}.  So
similar theorems can be churned out almost without 
human intervention, as we will see below.
\end{remark}

\section{Theorem 3.5 of \kaw}

Now we are ready to prove Theorem 3.5 of \kaw, namely:
\begin{theorem}
\leavevmode
\begin{enumerate}[(i)]
\item $L+U = \Enn - \{ F_n - 1 \suchthat n \geq 1 \}$;
\item $L_0 + U = \Enn - \{ F_n - 1 \suchthat n \geq 1 \text{ and } n 
\text{ odd } \}$;
\item $L+ U_0 = \Enn - \{ F_n - 1 \suchthat n \geq 1 \text{ and } n
\text{ even } \}$.
\end{enumerate}
\end{theorem}

\begin{proof}
We use the {\tt Walnut} commands
\begin{verbatim}
def thm35i "?msd_fib ~Ea,b (n=a+b) & $lower(a) & $upper(b)":
def thm35ii "?msd_fib ~Ea,b (n=a+b) & $lower0(a) & $upper(b)":
def thm35iii "?msd_fib ~Ea,b (n=a+b) & $lower(a) & $upper0(b)":
\end{verbatim}
which compute automata for the {\it complements} 
$\Enn - (L+U)$, $\Enn - (L_0 + U)$, and $\Enn - (L+U_0)$. 
These commands produce the automata depicted below.  
By comparison with Figure~\ref{fig1}, 
\kaw's Theorem 3.5 immediately follows.
\begin{figure}[H]
\begin{center}
\includegraphics[width=3.1in]{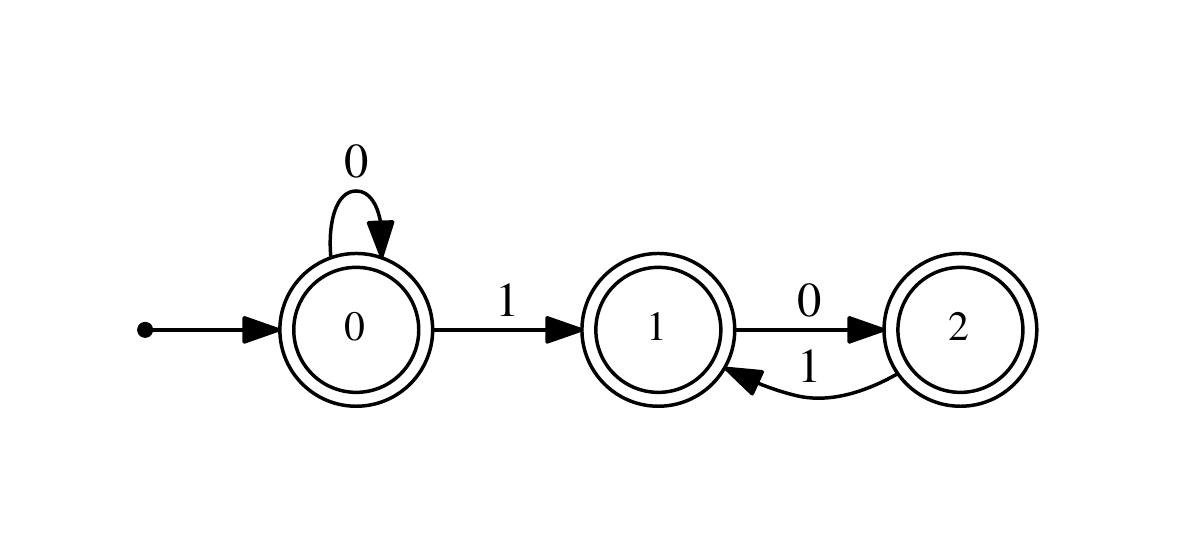}
\includegraphics[width=3.1in]{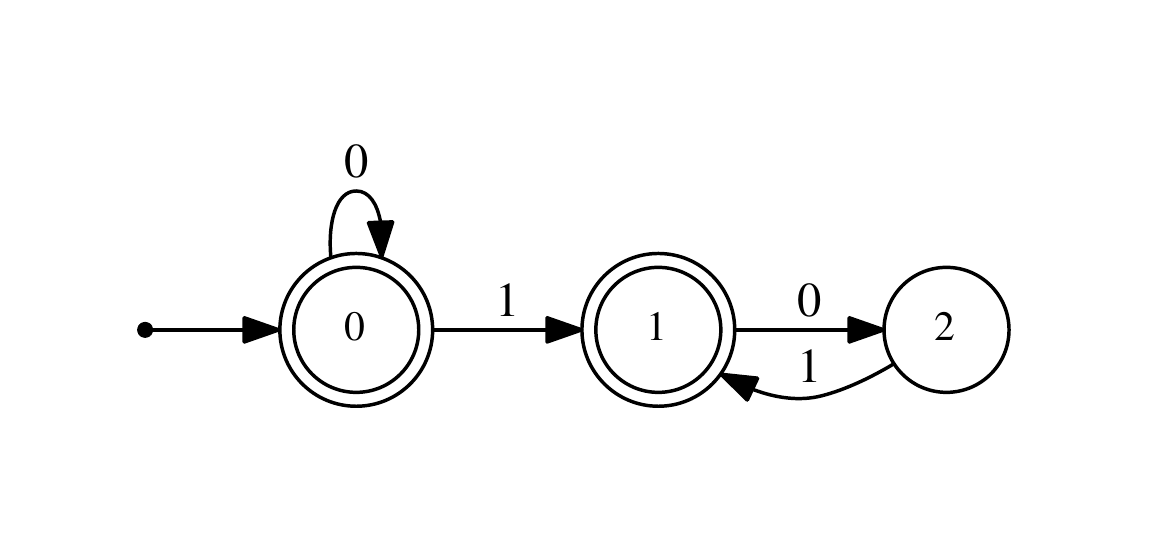}
\includegraphics[width=3.1in]{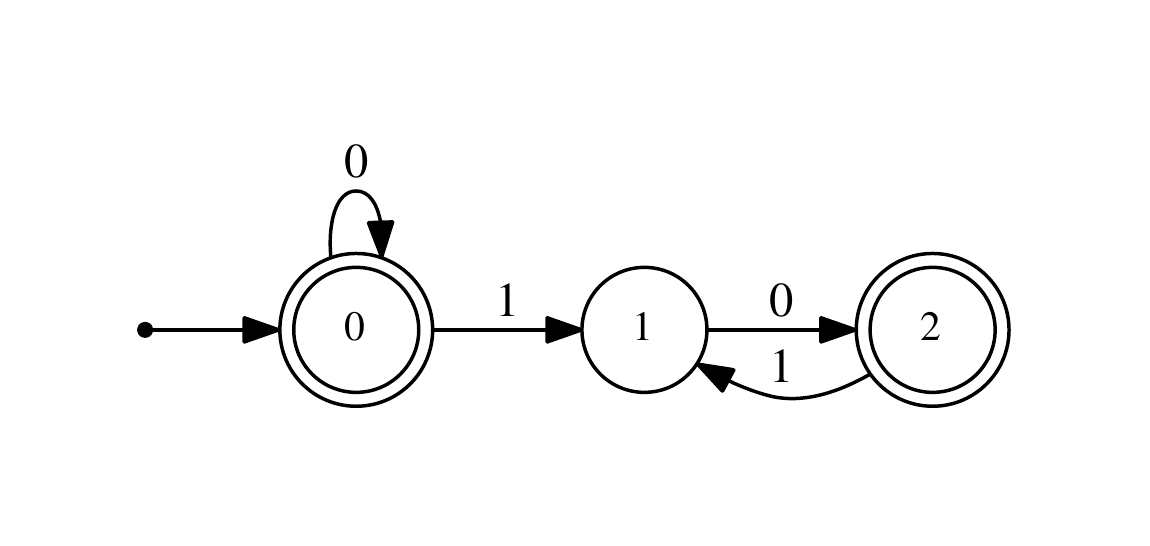}
\end{center}
\caption{Automata for $\Enn - (L+U)$, $\Enn - (L_0 + U)$, and $\Enn - (L+U_0)$}
\end{figure}
\end{proof}

\section{Theorem 3.7 of \kaw}

We now move on to Theorem 3.7 of \kaw:
\begin{theorem}
\leavevmode
\begin{enumerate}[(i)]
\item $L + U + U = \Enn - \{ 0,1,2,3,4,6,9 \}$;
\item $L_0 + U + U = \Enn - \{0,1,2,3,6\}$;
\item $L + U + U_0 = \Enn - \{0,1,2,4\}$;
\item $L + U_0 + U_0 = \Enn - \{0,2\}$.
\end{enumerate}
\label{thm37}
\end{theorem}

\begin{proof}
As above, we can easily create {\tt Walnut} predicates
for the complements of the three sumsets:
\begin{verbatim}
def thm37i "?msd_fib ~Ea,b,c (n=a+b+c) & $lower(a) & $upper(b) & $upper(c)":
def thm37ii "?msd_fib ~Ea,b,c (n=a+b+c) & $lower0(a) & $upper(b) & $upper(c)":
def thm37iii "?msd_fib ~Ea,b,c (n=a+b+c) & $lower(a) & $upper(b) & $upper0(c)":
def thm37iv "?msd_fib ~Ea,b,c (n=a+b+c) & $lower(a) & $upper0(b) & $upper0(c)":
\end{verbatim}
When we run these, we get the automata depicted below in
Figure~\ref{fig37}.  Since
each of these automata accept only a finite set, it is trivial to
verify the claims.
\begin{figure}[H]
\begin{center}
\includegraphics[width=3.1in]{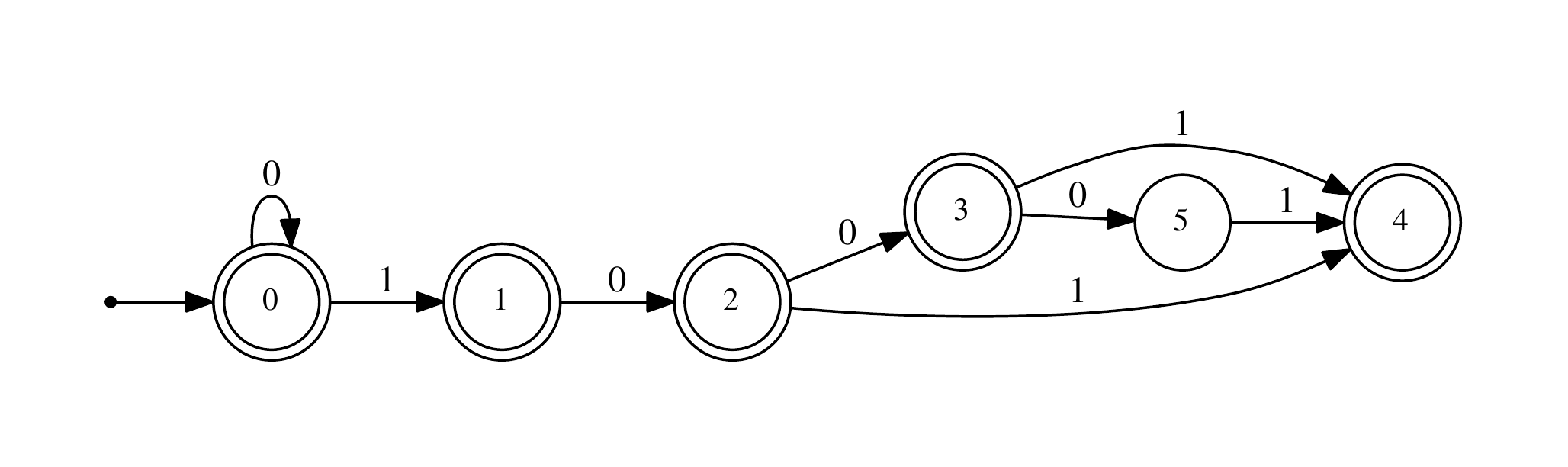}
\includegraphics[width=3.1in]{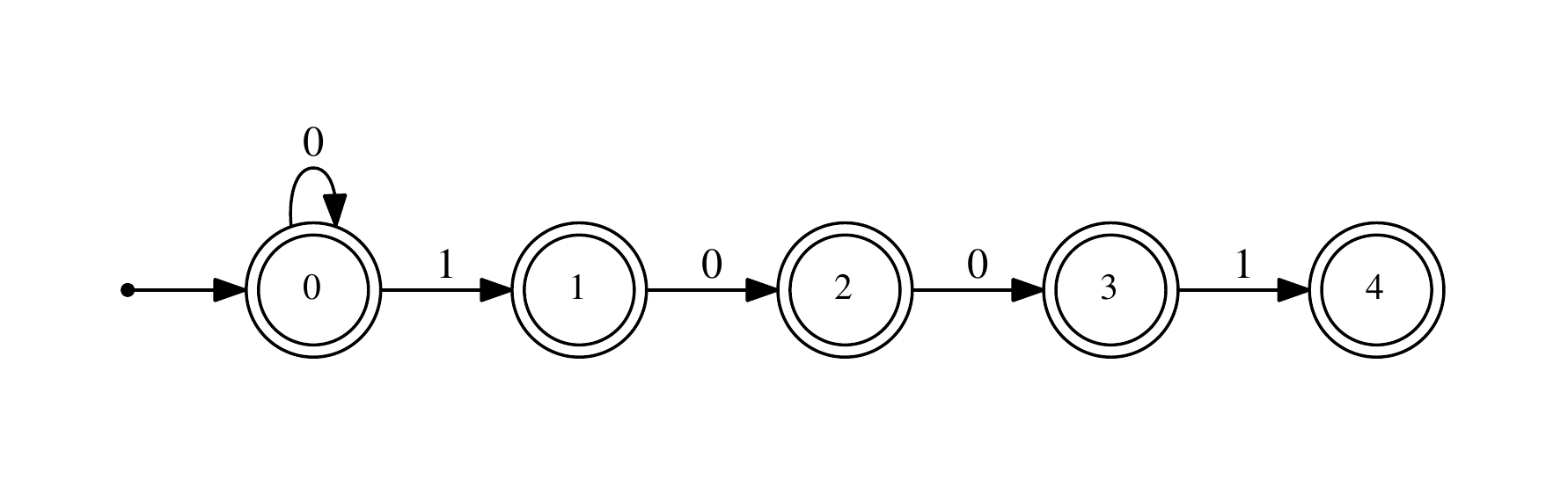}
\includegraphics[width=3.1in]{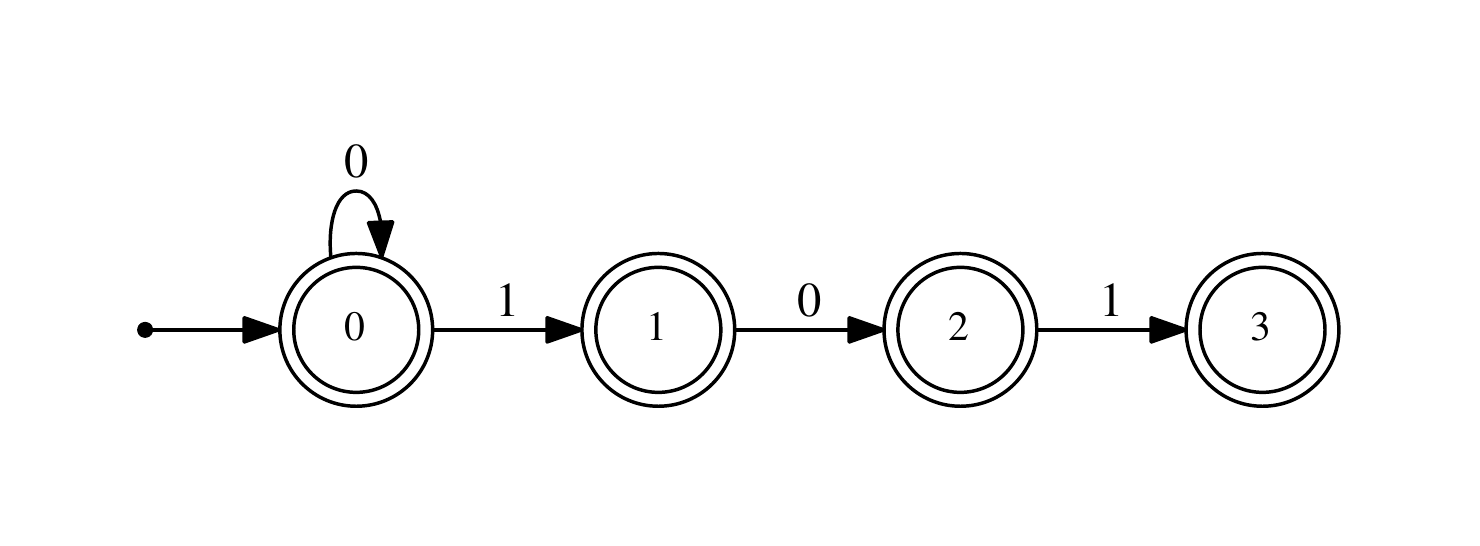}
\includegraphics[width=3.1in]{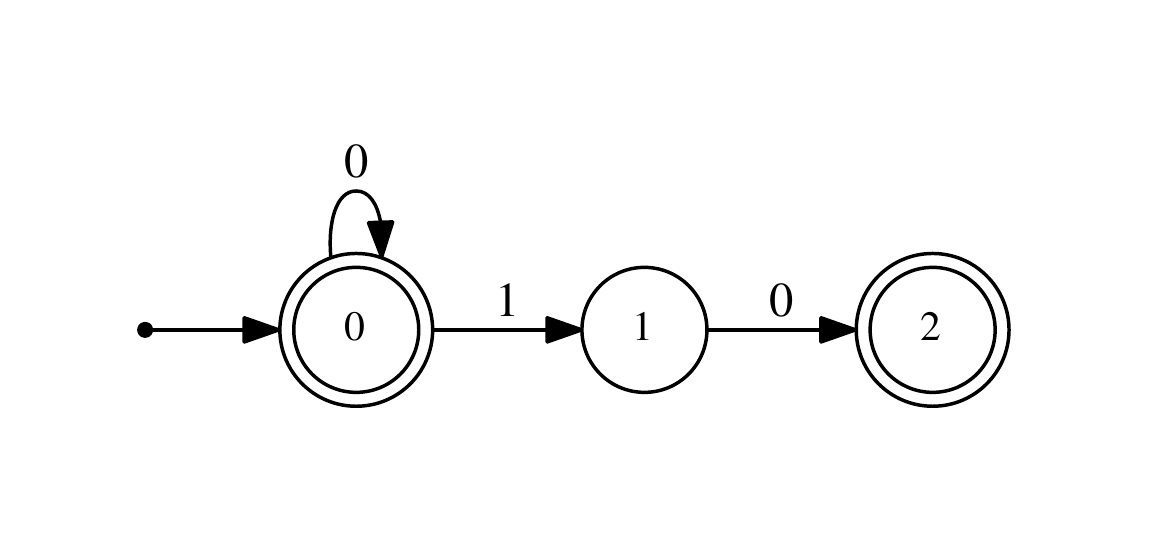}
\end{center}
\caption{Automata for parts (i)-(iv) of Theorem~\ref{thm37}}
\label{fig37}
\end{figure}
\end{proof}

\section{\kaw's Theorem 3.8}

\begin{theorem}
\leavevmode
\begin{enumerate}[(i)]
\item $U + U + U = \Enn - \{ 0,1,2,3,4, 5,7,8,10,13,18,26\} $;
\item $U_0 + U + U = \Enn - \{0,1,2,3,5,8,13\}$;
\item $U_0 + U_0 + U = \Enn - \{0,1,3,8\}$.
\end{enumerate}
\label{thm38}
\end{theorem}

\begin{proof}
As above, we can easily create {\tt Walnut} predicates
for the complements of the three sumsets:
\begin{verbatim}
def thm38i "?msd_fib ~Ea,b,c (n=a+b+c) & $upper(a) & $upper(b) & $upper(c)":
def thm38ii "?msd_fib ~Ea,b,c (n=a+b+c) & $upper0(a) & $upper(b) & $upper(c)":
def thm38iii "?msd_fib ~Ea,b,c (n=a+b+c) & $upper0(a) & $upper0(b) & $upper(c)":
\end{verbatim}
When evaluated in {\tt Walnut} these give the automata
below in Figure~\ref{fig38}.  Since
each of these automata accept only a finite set, it is trivial to
verify the claims.
\begin{figure}[H]
\begin{center}
\includegraphics[width=5.0in]{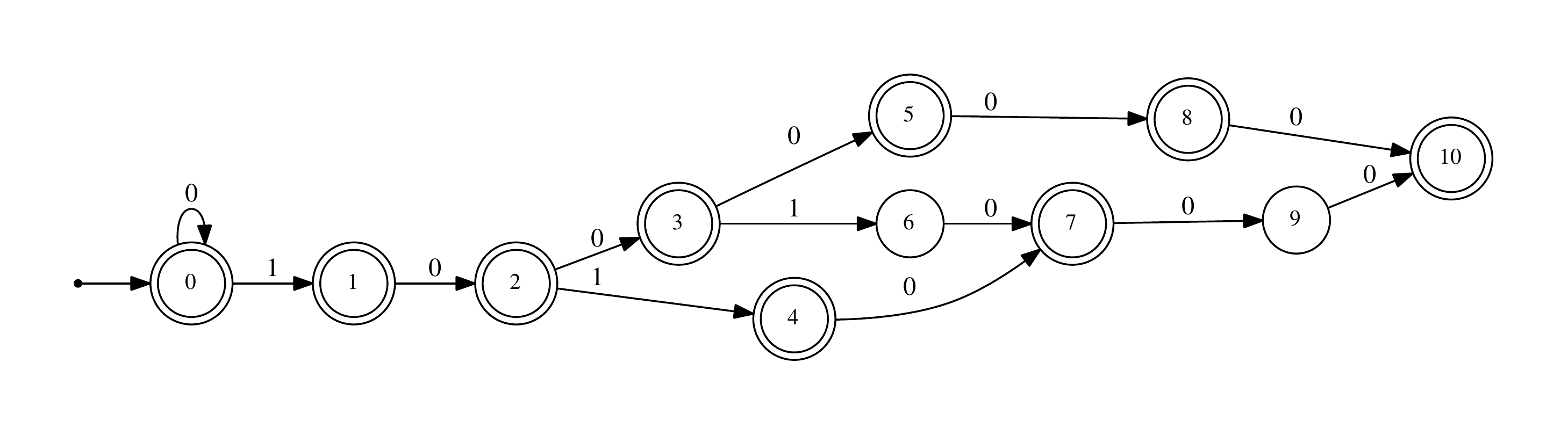}
\includegraphics[width=3.1in]{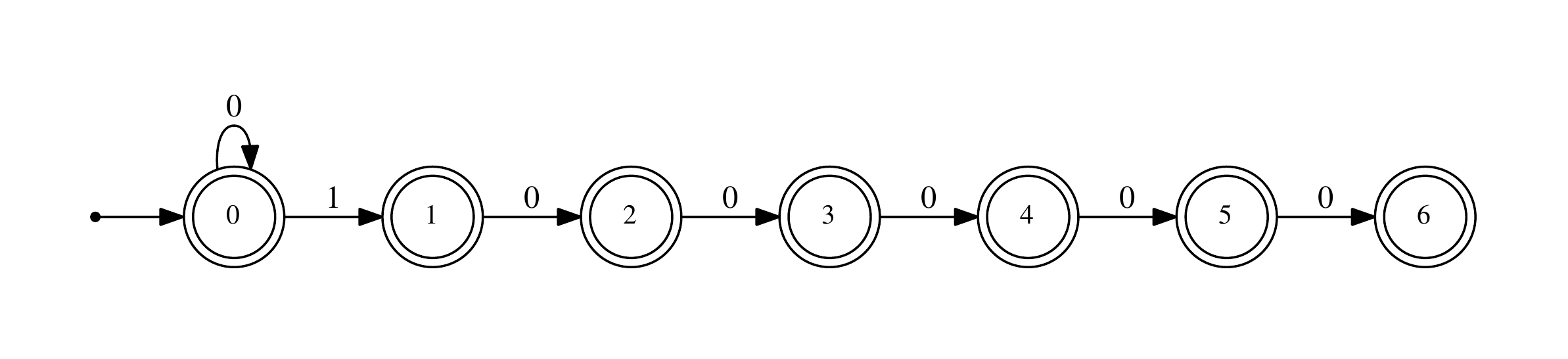}
\includegraphics[width=3.1in]{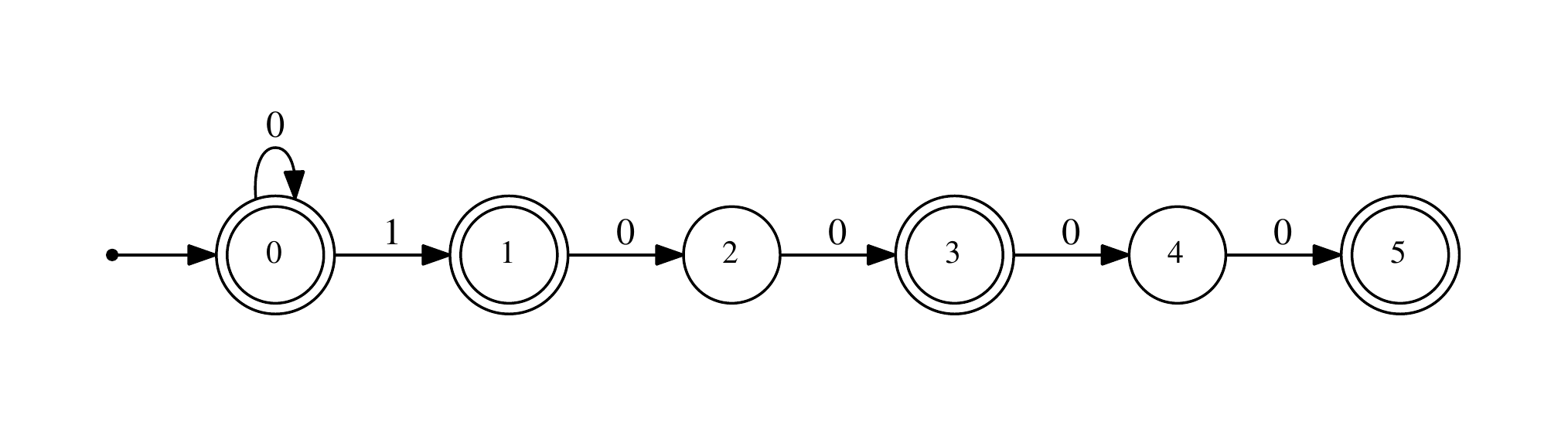}
\end{center}
\caption{Automata for parts (i)-(iii) of Theorem~\ref{thm38}}
\label{fig38}
\end{figure}
\end{proof}

\section{The sumset $U+U$}

Finally, perhaps the most interesting and
complicated object in the paper of \kaw~is
the sumset $U+U$.   They do not actually
provide a clean explicit definition
of this sumset, but with the method described here we can completely
describe it with an automaton
with 12 states, as follows:
\begin{theorem}
The sumset $U+U$ is accepted by the automaton in
Figure~\ref{figuu}.
\begin{figure}[H]
\begin{center}
\includegraphics[width=6.5in]{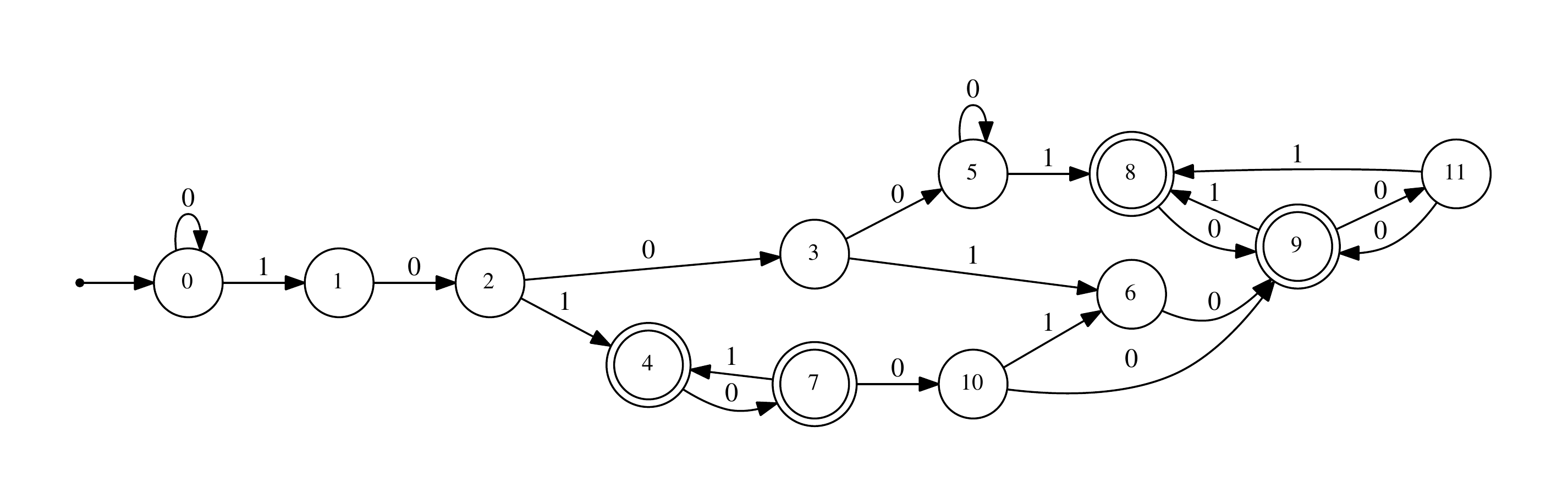}
\caption{Automaton for the sumset $U+U$}
\label{figuu}
\end{center}
\end{figure}
\end{theorem}

\begin{proof}
We use the {\tt Walnut} command:
\begin{verbatim}
def uplusu "?msd_fib Ea,b (n=a+b) & $upper(a) & $upper(b)":
\end{verbatim}
\end{proof}

With only a small amount of more work, from this
automaton all of the results in 
\kaw's Theorem 3.10 can be obtained.  We omit the details, as they
are just like what we have already done.

The facts that this automaton has 12 states and the sumset has an
infinite complement explains why it is hard to write down a
concise expression for $U+U$.  We argue that the automaton is
the best way to understand it.
Furthermore, the automaton also explains their remarks about the
``fractal'' properties of $U + U$, since automatic sets
are well-known to exhibit this kind of self-similarity 
\cite{Barbe&vonHaeseler:2003}.

\section{New results on the Wythoff sets}

\subsection{Counting representable numbers}

As an example of the power of the method, let's look at a problem that
\kaw~did not study:  how many numbers $<n$ are in the sumset
$U+U$?   Call this number $c(n)$.

We can easily create a
Fibonacci automaton that recognizes the pairs $(m,n)$ such that
$m < n$ and $m \in U+U$.
From this, the techniques described in
\cite{Du&Mousavi&Schaeffer&Shallit:2016} explain how to explicitly
compute a so-called ``linear representation'' for $c(n)$.  It can
be computed by the {\tt Walnut} command below:
\begin{verbatim}
eval count2upp n "?msd_fib (m < n) & E a,b (m=a+b) & $upper(a) & $upper(b)":
\end{verbatim}
In this case, the 
linear representation of rank $t$ that is
produced is a triple $(v, \mu, w)$, where
\begin{itemize}
\item $v$ is a $1 \times t$ matrix;
\item $\mu$ is a morphism from $\{ 0,1 \}$ to the set of $t \times t$
matrices with non-negative integer entries;
\item $w$ is a $t \times 1$ matrix,
\end{itemize}
such that $c(n) = v \mu( (n)_F) w$.
In this case, the rank is $33$; since it is so large, we do not display
it here.   However, it provides an efficient way to compute $c(n)$.

Furthermore, we can use it to explicitly compute $c(F_n)$:
\begin{theorem}
The number of elements of $U+U$ less than $F_n$ is $c(F_n)$ and
equals $2 F_{n-2} + 2 - n$ for $n \geq 4$.
\end{theorem}

\begin{proof}
We have $(F_n)_F = 1 0^{n-2}$ for $n \geq 2$.   Therefore
$c(F_n) = v \mu(1) \mu(0)^{n-2} w$ for $n \geq 2$.
If $M$ is a matrix, then it is well-known
that each entry of $M^i$ satisfies a linear recurrence
that is annihilated by the minimal polynomial of $M$.
Of course, the same thing
is also true for any linear combination of the entries
of $M^i$.  In this case, a symbolic algebra system such
as {\tt Maple} can easily compute the minimal polynomial of
$M = \mu(0)$; it is $X^3 (X+1)(X^2 - X - 1)(X-1)^2$.
Now the fundamental theorem of linear recurrences tells us
that $c(n)$ for $n \geq 5$ can be written as
$$ c_1 (-1)^n + c_2 \alpha^n + c_3 \beta^n + c_4 n^2 + c_5 n + c_6 $$
for some constants $c_1, c_2, \ldots, c_6$.
We can then determine the constants by evaluating $c(n)$ for
$5 \leq n \leq 10$ and solving the resulting linear system.
When we do this, and simplify, we get the desired result for $n \geq 5$.
We can then check that the formula also holds for $n = 4$.
\end{proof}

\subsection{Number of representations}

Given that a number is representable,
how many different representations are there
for $n$ as a sum of the form $a_i + a_j$, or $b_i + b_j$, or
$a_i + b_j$?
We can handle this problem using the enumeration capabilities of
{\tt Walnut}, as we did in the previous section.  We can compute
the corresponding linear representations using the commands below.
\begin{verbatim}
eval lowpluslow n "?msd_fib Eb (n=a+b) & $lower(a) & $lower(b)":
eval lowplusupp n "?msd_fib Eb (n=a+b) & $lower(a) & $upper(b)":
eval uppplusupp n "?msd_fib Eb (n=a+b) & $upper(a) & $upper(b)":
\end{verbatim}
These give us linear representations of rank $19$, from which
we can easily prove results like the following:
\begin{theorem}
The number of pairs $(i,j)$ such that 
$F_n - 1 = a_i + a_j$ is $F_{n-1} - 1$, for $n \geq 2$.
\end{theorem}

\section{Generalization to higher-order recurrences}

In their paper, \kaw~ask about the generalization of sumsets to
higher-order linear recurrences.   One such linear
recurrence is the famous
Tribonacci sequence \cite{Feinberg:1972}.
There is an analogue of the Fibonacci representation of integers,
called the Tribonacci representation, and having many of the
same properties.  We denote it as $(n)_T$.
Our automaton-based method also works for such sequences,
and Tribonacci-automatic sequences are built into {\tt Walnut}.
For more examples of {\tt Walnut} with the Tribonacci sequence,
see \cite{Mousavi&Shallit:2015}.

Here the appropriate analogues of the sequences $(a_n)_{n \geq 1}$ and
$(b_n)_{n \geq 1}$ we have considered so far in this paper
are the Carlitz-Scoville-Hoggatt ``sequences of higher
order'' $(A_n)_{n \geq 1}$, $(B_n)_{n \geq 1}$, and
$(C_n)_{n \geq 1}$; see \cite{Carlitz&Scoville&Hoggatt:1972,Dekking}.

To give some idea of what kinds of things can be proved,
consider the following two theorems:
\begin{theorem}
A number $n$ cannot be written in the form $A(i) + B(j)$ iff 
$n = 6$ or $(n)_T$ is a nonempty prefix of $100100100\cdots$.
\end{theorem}

\begin{theorem}
Every integer $n$, except 
$0,1,2,3,4,5,6,8,10,12,19$,
can be written in the form $n = A(i) + B(j) +C(k)$.
\end{theorem}

These, and many others, can be easily proved with {\tt Walnut}
almost immediately by mimicking the analysis we used above.
We give the {\tt Walnut} commands for these two theorems, and leave
further exploration as a fun exercise for the reader.
\begin{verbatim}
reg tend0 msd_trib "(0|1)*0":
reg tend01 msd_trib "1|((0|1)*01)":
reg tend011 msd_trib "11|((0|1)*011)":
def aa "?msd_trib Em $tend0(m) & n=m+1":
def bb "?msd_trib Em $tend01(m) & n=m+1":
def cc "?msd_trib Em $tend011(m) & n=m+1":
def aaplusbb "?msd_trib ~Ea,b (n=a+b) & $aa(a) & $bb(b)":
def aaplusbbpluscc "?msd_trib ~Ea,b,c (n=a+b+c) & $aa(a) & $bb(b) & $cc(c)":
\end{verbatim}

\section{Remarks}

The point we have tried to make by writing this paper is not to denigrate
the ingenuity or hard work of others.   It is, instead, to publicize
the fact that many 
case-based arguments in combinatorics on words in
the literature can now be replaced
by some simple computations, using a decision procedure
for the appropriate logical theory, using freely available
software.

The advantages to this approach are many:
\begin{itemize}
\item long ad hoc arguments can be handled in a simple and unified way;
\item drudgery is replaced by a software tool;
\item the software allows us to conduct experiments with ease,
giving us a ``telescope'' to view results that
at first appear only distantly provable;
\item in many cases {\tt Walnut}'s results actually provide the complete
statement of the desired result;
\item the formalism of finite
automata allow complicated results to be stated relatively simply;
\item when the automata are very complicated, likely there is no
really simple way to characterize the corresponding sets;
\item much more general results, such as enumeration, can also be obtained; and
finally,
\item generalizations to other settings, such as the
Tribonacci numbers here, are often easy and routine.
\end{itemize}
In this respect, the method described in this paper
can be viewed as another weapon
in the combinatorialist's arsenal, much like the Wilf-Zeilberger
method can be used for binomial coefficient identities and beyond
\cite{Petkovsek}.
We hope other authors will consider this approach in their work.

The free software {\tt Walnut} we have discussed is available 
at \\
\centerline{\url{https://github.com/hamousavi/Walnut}} \\
and a manual for its use is \cite{Mousavi:2016}.
The web page\\
\centerline{\url{https://cs.uwaterloo.ca/~shallit/walnut.html} } \\
on the author's home page contains references to many papers using
the approach we have described here.

\end{document}